\documentclass[final]{siamltex}

\usepackage{amsfonts}
\usepackage{amsmath}
\usepackage{amssymb}
\usepackage{color}
\usepackage{graphicx}
\usepackage{times}
\usepackage[round]{natbib}
\usepackage[T1]{fontenc}

\DeclareMathSymbol{\asymle}{\mathord}{AMSa}{"2E}

\newcommand {\Vector} {\mbox {vec}}
\newcommand {\sigmamin} {\sigma_{\min}}
\newcommand {\acoeff} {\alpha}
\newcommand {\bcoeff} {\beta}
\newcommand {\scaleA} {\mathcal {A}}
\newcommand {\scaleB} {\mathcal {B}}
\newcommand {\scaleP} {\mathcal {P}}
\newcommand {\scaleR} {\mathcal {R}}
\newcommand {\reals} {{\mathbb R}}

\newcommand {\vds} {\mbox {\bfseries \scshape v}}
\newcommand {\boldtheta} {{\boldsymbol \theta}}
\newcommand {\boldkappa} {{\boldsymbol \kappa^{}_{\boldsymbol 2}}}

\newcommand {\strutwidth} {0pt}

\newcommand {\col} {\mbox {\rm col}\hspace{0.05em}}
\newcommand {\DA} {\Delta A}

\newcommand {\DAx} {\Delta (Ax)}
\newcommand {\Db} {\Delta b}
\newcommand {\Dr} {\Delta r}
\newcommand {\Dx} {\Delta x}
\newcommand {\Dy} {\Delta y}
\newcommand {\Du} {\Delta u}
\newcommand {\Dv} {\Delta v}
\newcommand {\lessORequal} [1] {\mathop{\,\le\,}\limits^{#1}}

\newcommand {\setarray} {\relax}
\newcommand {\onetwo} [2] {\left[ \setarray \begin {array} {c c} #1& #2 \end {array} \right]}
\newcommand {\textonetwo} [2] {[ \begin {array} {c c} #1& #2 \end {array}]}
\newcommand {\twoone} [2] {\left[ \setarray \begin {array} {c} #1\\ #2 \end {array} \right]}
\newcommand {\twotwo} [4] {\left[ \setarray \begin {array} {c c} #1& #2\\ #3& #4 \end {array} \right]}
\newcommand {\JrA} {J_r [\Vector (A)]}
\newcommand {\JxA} {J_x [\Vector (A)]}

\title {Spectral Condition Numbers of Full Rank Linear Least Squares Residuals and Orthogonal Projections}

\title {Spectral Condition Numbers of Orthogonal Projections\\ and Full Rank Linear Least Squares Residuals\footnote {Please read and cite the corrected, published article that can be accessed through the DOI link on the arXiv page for this document.}}

\author{Joseph F. Grcar\thanks{6059 Castlebrook Drive, Castro Valley, CA 94552 USA (jfgrcar@comcast.net, or na.grcar@na-net.ornl.gov).}}

\date {}

\begin {document}

\maketitle

\begin {abstract} 
A simple formula is proved to be a tight estimate for the condition number of the full rank linear least squares residual with respect to the matrix of least squares coefficients and scaled $2$-norms. The tight estimate reveals that the condition number depends on three quantities, two of which can cause ill-conditioning. The numerical linear algebra literature presents several estimates of various instances of these condition numbers. All the prior values exceed the formula introduced here, sometimes by large factors.
\end {abstract}




\pagestyle{myheadings}
\thispagestyle{plain}
\markboth{JOSEPH F. GRCAR}{LINEAR LEAST SQUARES RESIDUAL}

\section {Introduction}

\subsection {Purpose}

Least squares residuals are quite important numerically. The residuals measure the quality of fits in regression analysis, and forming orthogonal projections is an essential step in many iterative algorithms for linear equations or matrix eigenvalues. 

This paper determines a tight estimate for the condition number of the residual in full rank least squares problems. Equivalently, the condition number of orthogonal projections into the span of linearly independent vectors is also estimated. The condition numbers are with respect to the matrix of least squares coefficients and with respect to scaled $2$-norms. The condition number of the residual, like the solution, is the value of an optimization problem that does not have an explicit formula but which does have a tight estimate. 

This introduction provides some background material. Section \ref {sec:definition} discusses the evaluation of condition numbers from Jacobian matrices. Section \ref {sec:description} describes the tight estimate of the condition number and provides an example; this material is appropriate for classroom presentation. Section \ref {sec:derivation} proves that the condition number varies from the estimate within a factor of $\sqrt 2$; the linear algebra is complicated but straightforward given an identity from a previous paper \citep {Grcar2009a}. Section \ref {sec:comparison} compares the results to the literature. Section \ref {sec:discuss} discusses the application to projections and to iterative algorithms. 

\subsection {Prior Work}

Conditioning with respect to perturbations of the matrix $A$ is the most interesting aspect of least squares problems,
\begin {equation}
\label {eqn:original-problem}
x = \arg \min_u \| b - A u \|_2 \qquad r = b - Ax \, .
\end {equation}
The condition numbers $\chi_x(A)$, of $x$ with respect to $A$ for various norms, have been studied in dozens of papers and books since \citet {GolubWilkinson1966} discovered the condition number for $2$-norms can depend on the square of the matrix condition number. Thus, it was equally surprising when \citet [p.\ 16, eqn.\ 7.7] {Bjorck1967a} discovered the conditioning of the residual is independent of the square.\footnote {Bj\"orck derived a bound for the sum of condition numbers with respect to $A$ and $b$, $\chi_r (A) + \chi_r (b)$.}  Even so, Bj\"orck's original formula turned out to be an overestimate. Roughly the same formula is still found in many textbooks (section \ref {sec:GVLH}).

\citet [p.\ 224, eqn.\ 5.4] {Wedin1973} derived a perturbation bound for the residual with respect to $A$ again for $2$-norms. This paper shows that Wedin's bound contains an estimate for $\chi_r (A)$ that is accurate within a factor of $2$. Wedin noted that his perturbation bound could be ``almost attained'' (p.\ 225). However, (p.\ 226) he also remarked that his paper only demonstrated near attainment for a perturbation bound on the least squares solution (not the residual). Thus the published literature has no prior proof of attainment for an error bound of the residual.

\citet {Geurts1982} and \citet {Gratton1996} have used Jacobian matrices to derive condition numbers, or estimates of condition numbers, for least squares solutions. Their results and those of \citet {Bjorck1967a}, \citet {Malyshev2003}, and \citet {Wedin1973} for the condition number of the solution are summarized by \citet {Grcar2009a}. There has been no similar determination of condition numbers based on Jacobian matrices for the residual. The spectral condition number of the residual, like the solution, is the value of an optimization problem that does not have an explicit formula but which does have a tight estimate. No tight estimates for the condition number of the least squares residual have been established previously. 

\section {Condition numbers}
\label {sec:definition}

\subsection {Error bounds and definitions of condition numbers}

``Perturbation bounds'' are used in numerical analysis to limit the sensitivity of the solution of a problem to changes in the initial data. Such bounds are customarily derived using matrix-vector algebra and norms; the coefficients of the data perturbations in these bounds are sometimes referred to as condition numbers.  For example, in one of the earliest books on rounding error analysis, \citet [p.\ 29] {Wilkinson1963} wrote ``we shall refer to [the coefficients] as condition numbers \dots.'' Many numerical analysts probably agree with Wilkinson in the interest of deriving error bounds, but the name ``condition number'' is used sparingly because the coefficients are only upper limits for condition numbers unless the error bounds are the smallest possible, equivalently, unless the error bounds are attained.   \citet [p.\ 1187] {Malyshev2003} observed, ``the bounds are commonly accepted as condition numbers, and any discussion about their sharpness is usually avoided.'' 

The oldest way to derive perturbation bounds is by differential calculus. If $y = f(x)$ is a vector valued function of the vector $x$ whose partial derivatives are continuous, then the partial derivatives give the best estimate of the change to $y$ for a given change to $x$
\begin {equation}
\label {eqn:approximation-1}
\Dy = f(x + \Dx) - f(x) \approx J_f (x) \, \Dx
\end {equation}
where $J_f (x)$ is the Jacobian matrix of the partial derivatives of $y$ with respect to $x$. The magnitude of the error in the first order approximation (\ref {eqn:approximation-1}) is bounded by Landau's little $o ( \| \Dx \| )$ for all sufficiently small $\| \Dx \|$.\footnote {The $o ( \| \Dx \| )$\ agreement is independent of the norm because all norms for finite dimensional spaces are equivalent \citep [p.\ 54, thm.\ 1.7] {Stewart1990}.} Thus $J_f (x) \, \Dx$ is the unique linear approximation to $\Dy$ in the vicinity of $x$.\footnote {Any other linear function added to $J_f (x) \, \Dx$ differs from $\Dy$ by ${\mathcal O} (\| \Dx \|)$ and therefore does not provide a $o ( \| \Dx \| )$ approximation.} Taking norms produces a perturbation bound,
\begin {equation}
\label {eqn:calculus-1}
\| \Dy \| \le \| J_f (x) \| \, \| \Dx \| + o (\| \Dx \|) \, .
\end {equation}
Equation (\ref {eqn:calculus-1}) is the smallest possible bound on $\| \Dy \|$ in terms of $\| \Dx \|$ provided the norm for the Jacobian matrix is induced from the norms for $\Dy$ and $\Dx$. In this case for each $x$ there is some $\Dx$, which is nonzero but may be chosen arbitrarily small, so the bound (\ref {eqn:calculus-1}) is attained to within the higher order term, $o (\| \Dx \|)$. There may be many other ways to define condition numbers, but because equation (\ref {eqn:calculus-1}) is the smallest possible bound, any definition of a condition number for use in bounds equivalent to (\ref {eqn:calculus-1}) must arrive at the same value, $\chi_y (x) = \| J_f (x) \|$.\footnote {A theory of condition numbers in terms of Jacobian matrices was developed by \citet [p.\ 292, thm.\ 4] {Rice1966}. See also \citet [p.\ 90] {Trefethen1997} for the present definition.} The matrix norm may be too complicated to have an explicit formula, but tight estimates can be derived as in this paper. 

\subsection {One or separate condition numbers}
\label {sec:separate}

Many problems depend on two parameters $u$, $v$ which may consist of the entries of a matrix and a vector (for example). In principle it is possible to treat the parameters altogether.\footnote {As will be discussed, \citet {Gratton1996} derived a joint condition number of the least squares solution with respect to a Frobenius norm of the matrix and vector that define the problem.}\footnote {\citet {Gratton1996} derived a joint condition number of the least squares solution with respect to a Frobenius norm of the matrix and vector that define the problem.} A condition number for $y$ with respect to joint changes in $u$ and $v$ requires a common norm for perturbations to both. Such a norm is 
\begin {equation}
\label {eqn:joint-norm}
\max \big\{ \| \Du \|, \, \| \Dv \| \big\} \, .
\end {equation}
A single condition number then follows that appears in an optimal error bound,
\begin {equation}
\label {eqn:single}
\| \Dy \| \le   \| J_f (u, v) \| \, \max \big\{ \| \Du \|, \, \| \Dv \| \big\} + o \left(\max \big\{ \| \Du \|, \, \| \Dv \| \big\} \right) .
\end {equation}
The Jacobian matrix $J_f (u, v)$ contains the partial derivatives of $y = f(u,v)$ with respect to the entries of both $u$ and $v$. The value of the condition number is again $\chi_y(u, v) = \| J_f (u, v) \|$ where the matrix norm is induced from the norm for $\Dy$ and the norm in equation (\ref {eqn:joint-norm}). 

Because $u$ and $v$ may enter into the problem in much different ways, it is customary to treat each separately.  This approach recognizes that the Jacobian matrix is a block matrix
\begin {displaymath}
J_f (u, v) = \onetwo {J_{f_1} (u)} {J_{f_2} (v)}
\end {displaymath}
where the functions $f_1 (u) = f(u, v)$ and $f_2(v) = f(u, v)$ have $v$ and $u$ fixed, respectively. 
The first order differential approximation (\ref {eqn:approximation-1}) is unchanged but is rewritten with separate terms for $u$ and $v$,
\begin {equation}
\label {eqn:approximation-2}
\Dy \approx J_{f_1} (u) \, \Du + J_{f_2} (v) \, \Dv \, ,
\end {equation}
and a perturbation bound is obtained by applying the triangle inequality,
\begin {eqnarray}
\nonumber
\| \Dy \|& \le& \| J_{f_1} (u)  \Du + J_{f_2} (v) \Dv \| + o \left(\max \big\{ \| \Du \|, \, \| \Dv \| \big\} \right)\\
\noalign {\smallskip}
\label {eqn:double}
& \le& \| J_{f_1} (u) \|  \, \| \Du \| + \| J_{f_2} (v) \| \, \| \Dv \| + o \left(\max \big\{ \| \Du \|, \, \| \Dv \| \big\} \right) \, .
\end {eqnarray}
The coefficients $\chi_y(u) = \| J_{f_1} (u) \|$ and $\chi_y(v) = \| J_{f_2} (v) \|$ are the separate condition numbers of $y$ with respect to $u$ and $v$, respectively. 

These two different approaches lead to error bounds (\ref {eqn:single}, \ref {eqn:double}) that differ by at most a factor of $2$. This fact is a property of induced norms. Consider a $p \times (m + n)$ block matrix
\begin {displaymath}
y = \onetwo A B \twoone u v
\end {displaymath}
and suppose norms are given for $\reals^p$, $\reals^m$ and $\reals^n$ as spaces of column vectors. A norm can be defined for $\reals^{m+n}$ as
\begin {displaymath}
\left\| \twoone u v \right\| = \max \big\{ \| u \|, \, \| v \| \big\} \, .
\end {displaymath}
These norms for $\reals^p$, $\reals^m$, $\reals^n$, and $\reals^{m+n}$ induce norms for $A$, $B$, and $\textonetwo A B$,
\begin {displaymath}
\| A \| = \max_{u \ne 0} {\| A u \| \over \| u \|} \, , \quad 
\| B \| = \max_{v \ne 0} {\| A u \| \over \| v \|} \, , \quad
\left\| \onetwo A B \right\| = \max_{u \ne 0 \; \mbox {\scriptsize or} \; v \ne 0}  {\| A u + B v \| \over \max \big\{ \| u \|, \, \| v \| \big\}} \, .
\end {displaymath}
The norm of the block matrix has a simple upper bound,
\begin {eqnarray}
\nonumber
\left\| \onetwo A B \right\|& =& \max_{u \ne 0 \; \mbox {\scriptsize or} \; v \ne 0}  {\| A u + B v \| \over \max \big\{ \| u \|, \, \| v \| \big\}}
\\ \noalign {\smallskip} 
\nonumber
& \le& \max_{u \ne 0 \; \mbox {\scriptsize or} \; v \ne 0}  {\| A u \| \over \max \big\{ \| u \|, \, \| v \| \big\}} + \max_{u \ne 0 \; \mbox {\scriptsize or} \; v \ne 0}  {\| B v \| \over \max \big\{ \| u \|, \, \| v \| \big\}}
\\ \noalign {\smallskip}
\nonumber
& =& \max_{u \ne 0}  {\| A u \| \over \| u \|} + \max_{v \ne 0}  {\| B v \| \over \| v \|}
\\ \noalign {\smallskip}
\label {eqn:lemma1}
& =& \| A \| + \| B \| \, ,
\end {eqnarray}
and a simple lower bound,
\begin {eqnarray}
\nonumber
\| A \| = \max_{u \ne 0}  {\| A u \| \over \| u \|}
& =& \max_{u \ne 0 \; \mbox {\scriptsize and} \; v = 0}  {\| A u + B v \| \over \max \big\{ \| u \|, \, \| v \| \big\}}
\\ \noalign {\smallskip}
\label {eqn:lemma2}
& \le& \max_{u \ne 0 \; \mbox {\scriptsize or} \; v \ne 0}  {\| A u + B v \| \over \max \big\{ \| u \|, \, \| v \| \big\}}
= \left\| \onetwo A B \right\| \, ,
\end {eqnarray}
and similarly $\| B \| \le \| \textonetwo A B \|$. Altogether, from equations (\ref {eqn:lemma1}, \ref {eqn:lemma2}),
\begin {equation}
\label {eqn:lemma3}
{\| A \| + \| B \| \over 2} \le \max \big\{ \| A \|, \, \| B \| \big\} \le \left\| \onetwo A B \right \| \le \| A \| + \| B \| \,
\end {equation}
which means that $\| A \| + \| B \|$ overestimates $\| \textonetwo A B \|$ by at most a factor of $2$. Returning to the Jacobian matrices $A = J_{f_1} (u)$, $B = J_{f_2} (v)$, and $\textonetwo A B = J_f (u, v)$, equation (\ref {eqn:lemma3}) can be rewritten
\begin {equation}
\label {eqn:sum-6}
{\chi_y (u) + \chi_y (v) \over 2} \le \chi_y (u, v) \le \chi_y (u) + \chi_y (v) \, .
\end {equation}
Thus, for the purpose of deriving tight estimates of joint condition numbers, it suffices to consider $\chi_y (u)$ and $\chi_y (v)$ separately. 

\section {Conditioning of the least squares residual}
\label {sec:description}

\subsection {Reason for considering full rank problems}
\label {sec:reason}

For any matrix $A$ and any similarly sized column vector $b$, the linear least squares problem (\ref {eqn:original-problem}) need not have an unique solution $x$, but it always has an unique residual $r = b - Ax = (I - P) b$ where $P = A A^\dag$ is the orthogonal projection into the column space of $A$, $\col (A)$, and where $A^\dag$ is the pseudoinverse of $A$. If $A$ has full column rank, then $P = A (A^t A)^{-1} A^t$. Changes to $A$ affect $r$ differently when $A$ does not have full column rank. In that case, $b \in \col (A + b z^t)$ for every nonzero right null vector $z$, so small changes to $A$ can produce large changes to $r$. In other words, a condition number of $r$ with respect to rank deficient $A$ does not exist or is ``infinite.'' Perturbation bounds in the rank deficient case can be found by restricting changes to the matrix, for which see \citet [p.\ 30, eqn.\ 1.4.27] {Bjorck1996} and \citet [pp.\ 136--162] {Stewart1990}. That theory is beyond the scope of the present discussion.

\subsection {The condition numbers}
\label {sec:brief}

This section summarizes the results and presents an example. Proofs are in section \ref {sec:derivation}. It is assumed that $A$ has full column rank and neither the solution $x$ nor the residual $r$ of the least squares problem are zero. The residual is proved to have a condition number $\chi_r(A)$ with respect to $A$ within the limits,
\begin {equation}
\label {eqn:chiA}
{1 \over \sqrt 2} \, \fbox {$\displaystyle {\boldkappa \sqrt { 1 + \left( \cot (\boldtheta) \over \vds \right)^2}}$} \; \le \; \chi_r (A) \; \le \; \fbox {$\displaystyle {\boldkappa \sqrt { 1 + \left( \cot (\boldtheta) \over \vds \right)^2}}$} \, .
\end {equation}
The quantities $\boldkappa$, $\boldtheta$, and $\vds$ are written bold to emphasize they are the only quantities affecting the tight estimate of the condition number; they are defined below.
There is also a condition number with respect to $b$,
\begin {equation}
\label {eqn:chib}
\chi_r (b) = \fbox {$\csc (\boldtheta)$} \, .
\end {equation}
These are condition numbers when the following scaled $2$-norms are used to measure the perturbations to $A$, $b$, and $x$, 
\begin {equation}
\label {eqn:specific-scaled-norms}
{\| \DA \|_2 \over \| A \|_2} \, ,
\qquad
{\| \Db \|_2 \over \| b \|_2} \, ,
\qquad
{\| \Dr \|_2 \over \| r \|_2} \, .
\end {equation}
Like equation (\ref {eqn:double}), the two condition numbers appear in error bounds of the form,\footnote {The constant denominators $\| A \|_2$ and $\| b \|_2$ could be discarded from the $o$ terms because only the order of magnitude of the terms is pertinent.} 
\begin {equation}
\label {eqn:error-bound}
{\| \Dr \|_2 \over \| r \|_2} \le \chi_r (A) {\| \DA \|_2 \over \| A \|_2} + \chi_r (b) {\| \Db \|_2 \over \| b \|_2} + o \left( \max \left\{ {\| \DA \|_2 \over \| A \|_2}, \, {\| \Db \|_2 \over \| b \|_2} \right\} \right) ,
\end {equation}
where $r + \Dr$ is the residual of the perturbed problem,
\begin {equation}
\label {eqn:perturbed-problem}
x + \Dx = \arg \min_u \| (b + \Db) - (A + \DA) u \|_2 \, .
\end {equation}

\begin {figure} 
\centering 
\includegraphics [scale=1] {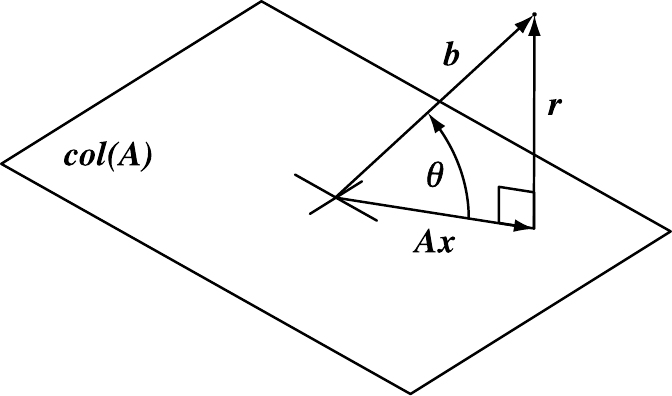}
\caption {Schematic of the least squares problem, the projection $Ax$,  and the angle $\boldtheta$ between $Ax$ and $b$.}
\label {fig:schematic}
\end {figure}

The quantities in the formulas (\ref {eqn:chiA}, \ref {eqn:chib}) are
\begin {equation}
\label {eqn:three}
\boldkappa = {\| A \|_2 \over \sigmamin} \, , \qquad 
\cot (\boldtheta) = {\|  A x \|_2 \over \| r \|_2} \, \qquad 
\vds = {\| A x \|_2 \over \| x \|_2 \, \sigmamin} \, ,
\end {equation}
where $\boldkappa$ is the spectral matrix condition number of $A$
($\sigmamin$ is the smallest singular value of $A$),
$\vds$ is van der Sluis's ratio between $1$ and $\boldkappa$,\footnote {\Citet [p.\ 251] {vanderSluis1975} introduced no notation. The Greek letters that look and sound like English v are $\nu$ and $\beta$, respectively, so it seems best to choose Roman $\vds$ for van der Sluis.} and $\boldtheta$ is the angle between $b$ and $\col (A)$. 
\smallskip
\begin {enumerate}
\item $\boldkappa$ depends only on the extreme singular values of $A$.
\item $\boldtheta$ depends only on the ``angle of attack'' of $b$ with respect to $\col (A)$.
\item If $A$ is fixed, then $\vds$ depends on the orientation of $b$ to $\col (A)$ but not on $\boldtheta$.\footnote {Because $A$ has full column rank, $Ax$ and $x$ can only vary proportionally when their directions are fixed.}
\end {enumerate}
\smallskip
Please refer to Figure \ref {fig:schematic} as needed. 
If $\col (A)$ is fixed, then $\boldkappa$ and $\boldtheta$ are separate sources of ill-conditioning for the residual. The ratio $\vds$ can never cause ill-conditioning because it only appears in the denominator of equation (\ref {eqn:chiA}) and $\vds$ is always at least $1$. Indeed, if $Ax$ has comparatively large components in singular vectors corresponding to the largest singular values, then $\vds \approx \boldkappa$ and $\vds$ might lessen the ill-conditioning caused by a small $\boldtheta$. 

\subsection {Conditioning example}
\label {sec:example}

This example illustrates the effects of $\boldkappa$, $\vds$, and $\boldtheta$ on $\chi_r (A)$. It is based on the example of \citet [p.\ 238] {Golub1996}.  Let
\begin {displaymath}
A = \left[ \begin {array} {c c} 1& 0\\ 0& \acoeff\\ 0& 0\\ \end {array}
\right] ,
\quad 
b = \left[ \begin {array} {c} \bcoeff \cos (\phi)\\ \bcoeff \sin (\phi) \\ 1 \end {array} \right] \, ,
\quad 
\DA = \left[ \begin {array} {r r} 0& 0\\ 0& 0\\ - \epsilon & - \epsilon \end {array} \right] ,
\end {displaymath}
where $0 < \epsilon \ll \acoeff, \bcoeff$, and $\acoeff < 1$. In this example, 
\begin {displaymath}
x = \left[ \begin {array} {c} \bcoeff \cos (\phi) \\ {\bcoeff \vphantom {(} \over \vphantom {(} \acoeff} \sin (\phi) \\ \end {array} \right] ,
\quad 
r = \left[ \begin {array} {c} 0\\ 0\\ 1 \end {array} \right] ,
\quad 
\Dr  = \left[ \begin {array} {c} 1 \\ {1 \over \acoeff} \\ \bcoeff \cos (\phi) + {\bcoeff \over \acoeff} \sin (\phi) \end {array} \right] \epsilon + {\mathcal O} (\epsilon^2)\, .
\end {displaymath} 
The three terms in the condition number are
\begin {displaymath}
\boldkappa = {1 \over \acoeff} \, , 
\qquad 
\vds = {1 \over \sqrt {[\acoeff \cos (\phi)]^2 + [\sin (\phi)]^2}} \, ,
\qquad
\cot (\boldtheta) = \bcoeff \, .
\end {displaymath}
These values can be manipulated by choosing $\acoeff$, $\bcoeff$ and $\phi$.  The tight upper bound on the condition number with respect to $A$ is
\begin {displaymath} 
\chi_r (A) \le {1 \over \acoeff} \, \sqrt {\strut 1 + [\acoeff \bcoeff \cos (\phi)]^2 + [\bcoeff \sin (\phi)]^2 } \, .
\end {displaymath}
The relative change to the residual,
\begin {displaymath}
{\| \Dr \|_2 \over \| r \|_2} 
= {1 \over \acoeff} \sqrt {\strut 1 + \acoeff^2 + [\acoeff \bcoeff \cos (\phi) + \bcoeff \sin (\phi)]^2} \, \epsilon + {\mathcal O} (\epsilon^2) \, ,
\end {displaymath}
can made be close to the bound on $\chi_r (A)$ times $\| \DA \|_2 / \| A \|_2 = \sqrt 2 \, \epsilon$. These formulas have been verified using Mathematica \citep {Wolfram2003}, as have formulas throughout the paper. 

\section {Derivation of the condition number estimates}
\label {sec:derivation}

\subsection {Choice of Norms}
\label {sec:norms}

\newcommand {\normA} [1] {\| #1 \|_{\scaleA}}
\newcommand {\normb} [1] {\| #1 \|_{\scaleB}}
\newcommand {\normr} [1] {\| #1 \|_{\scaleR}}

In theoretical numerical analysis especially for least squares problems the $2$-norm is preferred because for it the matrix condition number of $A^t A$ is the square of the matrix condition number of $A$. The norms used in this paper and in many other papers are defined as,
\begin {equation}
\label {eqn:norms}
\normA {\Vector (\DA)} = {\| \DA \|_2 \over \scaleA} \, ,
\qquad
\normb {\Db} = {\| \Db \|_2 \over \scaleB} \, ,
\qquad
\normr {\Dr} = {\| \Dr \|_2 \over \scaleR} \, ,
\end {equation}
where the choice of scale factors is left open. The scaling makes the size of the changes relative to the particular problem of interest. The scaling used in equations (\ref {eqn:chiA}--\ref {eqn:specific-scaled-norms}) is
\begin {equation}
\label {eqn:scale-factors}
\scaleA = \| A \|_2 \, ,
\qquad
\scaleB = \| b \|_2 \, ,
\qquad
\scaleR = \| r \|_2 \, .
\end {equation}
Some authors prefer to measure the residual relative to $b$ by choosing $\scaleR = \| b \|_2$. Other authors have no scaling, $\scaleA = \scaleB = \scaleR = 1$. All of these cases are accommodated by the notation in equation (\ref {eqn:norms}). The effect of the choice for $\scaleR$ is discussed in section \ref {sec:residual-scaling}.

\subsection {Notation}

The formula for the Jacobian matrix $J_r (b)$ of the residual $r = [I - A (A^t A)^{-1} A^t] b$ with respect to $b$ is clear.\footnote {The notation of section \ref {sec:separate} would introduce a name, $f_2$, for the function by which $r$ varies with $b$ when $A$ is held fixed, $r = f_2 (b)$, so that the notation for the Jacobian matrix is then $J_{f_2} (b)$. This pedantry will be discarded here to write $J_r(b)$ for the matrix of partial derivatives of $r$ with respect to $b$ with $A$ held fixed.} For derivatives with respect to the entries of $A$, it is necessary to use the ``$\Vector$'' construction to order the matrix entries into a column vector; $\Vector (B)$ is the column of entries $B_{i,j}$ with $i,j$ in co-lexicographic order.\footnote {The alternative to placing the entries of matrices into column vectors is to use more general linear spaces and the Fr\'echet derivative. That approach seems unnecessarily abstract because the spaces have finite dimension.} The first order approximation (\ref {eqn:approximation-2}) is then
\begin {equation}
\label {eqn:total-differential}
\Dr = \JrA \, \Vector (\DA) + J_r (b) \, \Db + \mbox {higher order terms in $\DA$ and $\Db$}
\end {equation}
and upon taking norms
\begin {eqnarray}
\nonumber
\normr {\Dr}& \le& \normr {\JrA \, \Vector (\DA)} + \normr {J_r (b) \, \Db} + o \, ( \dots )
\\ \noalign {\medskip}
\label {eqn:differential-bound}
& \le& \underbrace {\| J_r [\Vector (\DA)] \|}_{\displaystyle \chi_r (A)} \, \normA {\DA} + \underbrace {\displaystyle \| J_r (b) \|}_{\displaystyle \chi_r (b)} \, \normb {\Db} + o \, ( \dots )
\end {eqnarray}
where the norms of the two Jacobian matrices are induced from $\normr {\cdot}$, $\normA {\cdot}$ and from $\normr {\cdot}$, $\normb {\cdot}$, respectively. The high order term in equation (\ref {eqn:differential-bound}) is $o ( \max \{ \normA {\DA}, \, \normb {\Db} \} )$ because from equation (\ref {eqn:calculus-1}) the norm $\max \{ \normA {\cdot}, \, \normb {\cdot} \}$ has been given to the space that jointly consists of matrices $A$ and vectors $b$.

\subsection {Condition number of {\itshape r\/} with respect to {\itshape b\/}} 
\label {sec:conditionwrtb}

For the orthogonal projection $P$ defined in section \ref {sec:reason}, from $r = (I - P) b$ follows $J_r (b) = I - P$, hence
\begin {equation}
\label {eqn:conditionwrtb}
\| J_r (b) \| 
= \max_{\Db} {\normr {J_r (b) \, \Db} \over \normb {\Db}} 
= \max_{\Db} {\displaystyle \left( {\| (I - P)  \Db \|_2 \over \scaleR} \right) \over \displaystyle \left( {\| \Db \|_2 \over \scaleB} \right)}
= {\scaleB \over \scaleR} \, ,
\end {equation}
which is equation (\ref {eqn:chib}) for the choice of scale factors in equation (\ref {eqn:scale-factors}).

\subsection {Condition number of {\itshape r\/} with respect to {\itshape A\/}} 
\label {sec:begin}

\newcommand {\JrxA} {\left[ 
\setlength {\arraycolsep} {0.25em} \begin {array} {ccccc}
x_1 I& x_2 I& \cdots& x_n I\\
e_1 r^t& e_2 r^t& \cdots & e_n r^t
\end {array} \right]}

Evaluating the condition number of the residual requires a formula for the Jacobian matrix $\JrA$. 
Differentiating the entries of 
\begin {displaymath}
r = \big[ I - A \left( A^t A \right)^{-1} \hspace {-0.33em} A^t \big] b
\end {displaymath}
by those of $A$ seems to be a daunting task. Instead, $\JrA$ is constructed from the total differential of the identity,
\begin {displaymath}
\twotwo I A {A^t} {0} \twoone r x - \twoone b {0} =  0 \, .
\end {displaymath}
Assuming $b$ is fixed because it already has been treated in section \ref {sec:conditionwrtb}, the total differential is
\begin {displaymath}
\twotwo I A {A^t} {0} \twoone {dr} {dx} + \JrxA \Vector (dA) = 0 \, ,
\end {displaymath}
where $x_i$ is the $i$-th entry of $x$ and where $e_i$ is the $i$-th column of the $n \times n$ identity matrix. Hence 
\begin {eqnarray}
\nonumber
\twoone {dr} {dx}& =& - {\twotwo I A {A^t} {0}}^{-1} \JrxA \Vector (dA) 
\\ \noalign {\smallskip}
\nonumber
& =& - \twotwo {I - P\;} {A (A^t A)^{-1}} {(A^t A)^{-1} A^t\;} {- (A^t A)^{-1}} \JrxA \Vector (dA) 
\\ \noalign {\smallskip}
\label {eqn:Jrx(A)}
& = & \twoone {\JrA} {\JxA} \Vector (dA) 
\end {eqnarray}
in which $P = A (A^t A)^{-1} A^t$ is the orthogonal projection into the column space of $A$. The two matrix blocks in equation (\ref {eqn:Jrx(A)}) are the Jacobian matrices of $r$ and $x$ as functions of the entries of $A$ with $b$ held fixed.

\subsection {Transpose formula for condition numbers}
\label {sec:transpose}

The desired condition number is the norm induced from the norms in equation (\ref {eqn:norms}). 
\begin {displaymath}
\setlength {\arraycolsep} {0.25em}
\begin {array} {r c l}
\| \JrA \|
& =
& \displaystyle \max_{\DA} {\normr {\JrA \, \Vector (\DA)} \over \normA {\Vector (\DA)}}
\\ \noalign {\bigskip}
& =
& \displaystyle {\scaleA \over \scaleR} \max_{\DA} {\| \JrA \, \Vector (\DA) \|_2 \over \| \DA \|_2} 
\end {array}
\end {displaymath}
The numerator and denominator are vector and matrix $2$-norms, respectively. If $A$ is an $m \times n$ matrix, then this maximization is a large problem with $mn$ degrees of freedom. The identity for the norm of the transposed operator can be applied to reduce the degrees to $m$, 
\begin {equation}
\label {eqn:induced}
\| \JrA \| =  {\scaleA \over \scaleR} \max_{\Dr} {\| \JrA^t \Dr \|_2^* \over \| \Dr \|_2^*} \, .
\end {equation}
Here, the identical norm for the transposed Jacobian matrix is induced from the duals of the $2$-norms for matrices and vectors. The vector $2$-norm is its own dual. The dual of the matrix $2$-norm is determined in \cite {Grcar2009a} to be the sum of the singular values of the matrix, including multiplicities. This norm is sometimes called the nuclear norm or the trace norm.

\subsection {Condition number of {\itshape r\/} with respect to {\itshape A\/}, continued} 
\label {sec:conditionwrtAcontinued}

The application of equation (\ref {eqn:induced}) requires the evaluation of the matrix-vector product in the numerator. Note that for any vectors $r^\prime$ and $x^\prime$,
\begin {displaymath}
\JrxA^t \twoone {r^\prime} {x^\prime} = \Vector [r^\prime x^t + r \, (x^\prime)^t] \, .
\end {displaymath}
With this identity it is now possible to compute, from equation (\ref {eqn:Jrx(A)}), 
\begin {eqnarray*}
\{ \JrA \}^t \Dr
& =
& \twoone {\JrA} {\JxA}^t \twoone {\Dr} {{0}}
\\
& =
& - \JrxA^t \twotwo {I - P\;} {A (A^t A)^{-1}} {(A^t A)^{-1} A^t\;} {- (A^t A)^{-1}}  \twoone {\Dr} {{0}}
\\
\noalign {\medskip}
& =
& - \JrxA^t  \twoone {(I - P) \Dr} {(A^t A)^{-1} A^t \, \Dr}
\\
\noalign {\medskip}
& =
& \Vector (u_1^{} v_1^t + u_2^{} v_2^t) \, ,
\end {eqnarray*}
in which
\begin {equation}
\label {eqn:u-and-v}
\setlength {\arraycolsep} {0.125em}
\begin {array} {r c l r c l}
u_1& =& (I-P) \Dr \qquad& 
v_1& =& x\\ \noalign {\smallskip}
u_2& =& r&
v_2& =& (A^t A)^{-1} A^t \, \Dr \, .
\end {array}
\end {equation}
Thus equation (\ref {eqn:induced}) is the following optimization,
\begin {eqnarray}
\nonumber
\setlength {\arraycolsep} {0.25em}
\| \JrA \|
& =
& \displaystyle {\scaleA \over \scaleR} \max_{\Dr} {\| u_1^{} v_1^t + u_2^{} v_2^t \|^*_2 \over \| \Dr \|_2}\, .
\\ \noalign {\medskip}
\label {eqn:induced-transpose}
& =
& \displaystyle {\scaleA \over \scaleR} \max_{\| \Dr \|_2 = 1} \| u_1^{} v_1^t + u_2^{} v_2^t \|^*_2 \, .
\end {eqnarray}

For ease of notation, let $g(\Dr)$ be the objective function in equation (\ref {eqn:induced-transpose}). In \cite {Grcar2009a} it is shown that
\begin {equation}
\label {eqn:numerator}
\begin {array} {l}
g(\Dr) = {}\\ \noalign {\smallskip}
\quad \sqrt {\strut \| u_1 \|_2^2 \, \| v_1 \|_2^2 + \| u_2 \|_2^2 \, \| v_2 \|_2^2 + 2 \, \| u_1 \|_2 \, \| v_1 \|_2 \, \| u_2 \|_2 \, \| v_2 \|_2 \, \cos (\boldtheta_u - \boldtheta_v)} \, ,
\end {array}
\end {equation}
where $\boldtheta_u$ is the angle between $u_1$ and $u_2$, and $\boldtheta_v$ is the angle between $v_1$ and $v_2$, and both angles should be taken from $0$ to $\pi$. Evaluating the maximum has two parts. 

The first step shows $\Dr$ can be restricted so that $\cos (\boldtheta_u - \boldtheta_v) \ge 0$. The vector $\Dr$ always could be decomposed into a component in $\col (A)$ and a component orthogonal to this subspace. Let the component inside $\col (A)$ be $a^\prime$. Further, the component outside can be decomposed into components parallel to $r$ and orthogonal to $r$, say $\gamma  r + r^\prime$ for some coefficient $\gamma$. With these choices to express $\Dr$,
\begin {displaymath}
\Dr = \gamma r + r^\prime + a^\prime \qquad \mbox {where} \quad a^\prime \in \col (A), \quad r^\prime \perp \col (A), \quad r^\prime \perp r
\end {displaymath}
the vectors in equation (\ref {eqn:u-and-v}) are
\begin {displaymath}
\setlength {\arraycolsep} {0.125em}
\begin {array} {r c l r c l}
u_1& =& \gamma r + r^\prime \qquad& 
v_1& =& x\\ \noalign {\smallskip}
u_2& =& r&
v_2& =& (A^t A)^{-1} A^t a^\prime
\end {array}
\end {displaymath}
and the angles are
\begin {eqnarray*}
\cos (\boldtheta_u)
& =
& {u_1^t u_2 \over \| u_1 \|_2 \, \| u_2 \|_2} = {\gamma \| r \|_2 \over \sqrt {\gamma^2 \| r \|_2^2 + \| r^\prime \|_2^2}}
\\ \noalign {\medskip}
\cos (\boldtheta_v)
& =
& {v_1^t v_2 \over \| v_1 \|_2 \, \| v_2 \|_2} = {x^t  (A^t A)^{-1} A^t a^\prime \over \| x \|_2 \, \| (A^t A)^{-1} A^t a^\prime \|_2}
\end {eqnarray*}
Thus the sign of $\gamma$ affects only the angle $\boldtheta_u$ in equation (\ref {eqn:numerator}), so it can be chosen to place $\boldtheta_u$ in the same quadrant as $\boldtheta_v$ (either from $0$ to $\pi / 2$, or from $\pi / 2$ to $\pi$) and hence $\cos (\boldtheta_u - \boldtheta_v) \ge 0$. This means the maximum of equation (\ref {eqn:induced}) can be restricted to those $\Dr$ for which 
\begin {equation}
\label {eqn:lowerandupper}
\begin {array} {l c c c r}
L(\Dr) = {}\\
\noalign {\smallskip}
\quad \sqrt {\strut \| u_1 \|_2^2 \, \| v_1 \|_2^2 + \| u_2 \|_2^2 \, \| v_2 \|_2^2}
& \le& g(\Dr)& \le&
\| u_1 \|_2 \, \| v_1 \|_2 + \| u_2 \|_2 \, \| v_2 \|_2 \quad\\
\noalign {\smallskip}
&&&& {} = U (\Dr) \, .
\end {array}
\end {equation}

The second step chooses $\Dr$ to maximize the upper bound $U(\Dr)$.  As before, the vector $\Dr$ always can be decomposed into a component in $\col (A)$ and a component in the orthogonal complement. Without loss of generality, assume $\Dr = \cos (\phi) r^{\prime \prime} + \sin (\phi) a^{\prime \prime}$ where $a^{\prime \prime}$ and $r^{\prime \prime}$ are unit vectors in $\col (A)$ and the complement, respectively, and where the coefficients are determined by an angle $\phi$ between $0$ and $\pi / 2$.\footnote {The coefficients $\cos (\phi)$ and $\sin (\phi)$ are non-negative so the choice of $\phi$ does not affect the choice of sign needed for equation (\ref {eqn:lowerandupper}).} The vectors in equation (\ref {eqn:u-and-v}) for this representation of $\Dr$ are
\begin {displaymath}
\setlength {\arraycolsep} {0.125em}
\begin {array} {r c l r c l}
u_1& =& \cos (\phi) r^{\prime \prime} \qquad& 
v_1& =& x\\ \noalign {\smallskip}
u_2& =& r&
v_2& =& \sin (\phi) (A^t A)^{-1} A^t a^{\prime \prime} \, .
\end {array}
\end {displaymath}
The largest $\| v_2 \|_2$ occurs when $a^{\prime \prime}$ is a left singular vector for the smallest singular value of $A$, $\sigmamin$, in which case $v_2 = (\sin (\phi) / \sigmamin) \, a^{\prime \prime}$; altogether
\begin {displaymath}
U (\Dr) = \cos (\phi) \, \| x \|_2 + \sin (\phi) \, {\| r \|_2 \over \sigmamin} \, .
\end {displaymath}
The maximum of this formula with respect to $\phi$ determines an optimal $\Dr_{\rm bnd}$ where the upper bound is
\begin {equation}
\label {eqn:upper-bound}
U (\Dr_{\rm bnd}) = \sqrt {\left({\| r \|_2 \over \sigmamin}\right)^2 + \| x \|_2^2} \; .
\end {equation}
The maximum has been verified using Mathematica \citep {Wolfram2003}. 

The formula in equation (\ref {eqn:upper-bound}) is the maximum of the upper bounds, which is not to say it is the maximum of equation (\ref {eqn:induced}). The objective function $g$ and the lower and upper bounds $L$ and $U$, when evaluated at $\Dr_{\rm bnd}$ and $\Dr_{\max}$, must be arranged as follows,
\begin {displaymath}
L (\Dr_{\rm bnd})
\lessORequal a
g (\Dr_{\rm bnd}) 
\lessORequal b
g (\Dr_{\max}) 
\lessORequal c
U (\Dr_{\max})
\lessORequal d
U (\Dr_{\rm bnd}) \, .
\end {displaymath}
These inequalities have the following justifications:
(a) equation (\ref {eqn:lowerandupper}),
(b) choice of $\Dr_{\rm max}$,
(c) equation (\ref {eqn:upper-bound}), and
(d) choice of $\Dr_{\rm bnd}$.
Therefore equation (\ref {eqn:upper-bound}) is an upper bound for the maximum. From the formula for $L(\Dr)$ in equation (\ref {eqn:lowerandupper}), the upper bound is at most $\sqrt 2$ times larger than a lower bound for the maximum. Note that to complete the limits and the condition number, these values must be scaled by the coefficient $\scaleA / \scaleR$ in equation (\ref {eqn:induced}).

\subsection {Summary of condition numbers} 
\label {sec:conditionwrtAfinished}

\ \par

\begin {theorem} 
[\scshape Spectral condition numbers]
\label {thm:condition-numbers}
For the full rank linear least squares problem with solution $x = (A^tA)^{-1} A^t b$ and residual $r = b - Ax$, and for the scaled norms in equation (\ref {eqn:norms}) with scale factors $\scaleA$, $\scaleB$, and $\scaleR$,
\begin {equation}
\label {eqn:thm}
\setlength {\arraycolsep} {0.33em}
\begin {array} {r c c c l}
&& \displaystyle \chi_r (b)
& =
& \displaystyle {\scaleB \over \scaleR} \, ,
\\ \noalign {\bigskip}
\displaystyle {{1 \over \sqrt 2} \, {\scaleA \over \scaleR} \,\sqrt {\left({\| r \|_2 \over \sigmamin}\right)^2 + \| x \|_2^2}}
& \le 
& \displaystyle \chi_r (A)
& \le 
& \displaystyle {{\scaleA \over \scaleR} \,\sqrt {\left({\| r \|_2 \over \sigmamin}\right)^2 + \| x \|_2^2}} \, ,
\end {array}
\end {equation}
where $\sigmamin$ is the smallest singular value of $A$. These formulas simplify to those in section \ref {sec:brief} for the choice of scale factors in equation (\ref {eqn:scale-factors}).
\end {theorem}

\begin {proof} 
Section \ref {sec:conditionwrtb} derives $\chi_r(b)$, and sections \ref {sec:begin}--\ref {sec:conditionwrtAcontinued} derive the bounds on $\chi_r(A)$. 
\end {proof}

\section {Comparison with published bounds}
\label {sec:comparison}

\subsection {The estimate of Wedin} Table \ref {tab:published} lists the condition number estimates in some textbook error bounds for the least squares residual. All the values exceed the upper estimate of theorem \ref {thm:condition-numbers} to varying degrees. 

The very early formula of \citet [p.\ 224, eqn 5.4] {Wedin1973} is also reported in the more recent textbook of \citet [p.\ 30, eqn.\ 1.4.27] {Bjorck1996}. It is for perturbations only to $A$, that is for the choice $\Db = 0$, and for the choice of scale factors $\scaleA = \scaleR = 1$. The value exceeds the estimate in theorem \ref {thm:condition-numbers} by at most the factor $\sqrt 2$, so it is at most double the condition number. The other two values in Table \ref {tab:published} can be severe overestimates.

\begin {table} 
\caption {\it Condition number estimates in textbook error bounds for the least squares residual. The full rank least squares problem is $\min_x \| b - A x \|_2$, the solution is $x$, the residual is $r$, the smallest nonzero singular value of $A$ is $\sigma_{\min}$, the condition number of $A$ is $\boldkappa = \| A \|_2 / \sigmamin$.}
\label {tab:published}
\newcommand {\equationstrut} {\vrule depth1.0ex height2.5ex width\strutwidth}
\newcommand {\textstrut} {\vrule depth1.125ex height2.375ex width\strutwidth}
\newcommand {\leftbox} [1] {\begin {minipage}{8.5em}\footnotesize \vrule depth0ex height3.0ex width\strutwidth \raggedright #1\vrule depth1.5ex height2ex width\strutwidth \end {minipage}}
\newcommand {\twolines} [2] {\begin {array} {c} \displaystyle {#1}\\ \noalign {\smallskip} \hspace* {-1em} \displaystyle #2 \hspace*{-1em}\end {array}}
\newcommand {\fracstrut} {\vrule depth2ex height3.5ex width\strutwidth}
\setlength {\arraycolsep} {0.45em}
\vspace {-2ex}
\begin {displaymath}
\small
\begin {array} {| l | c | c | c | c | c |}
\cline {2-6}
\multicolumn {1} { c |} {\textstrut}
& \multicolumn {3} { c |} {\mbox {norms and scale factors}}
& \multicolumn {2} { r |} {\mbox {maximum overestimation factor}}
\\
\cline {1-5} 
\multicolumn {1} {| c |} {\textstrut \mbox {source}}
& \multicolumn {2} { c |} {\mbox {data}}
& \mbox {residual}
& \mbox {estimate for $\chi_r (A)^{**}$}
&
\\ \hline \hline
\leftbox {\citet [p.\ 224, eqn.\ 5.4] {Wedin1973}, \citet [p.\ 30, eqn.\ 1.4.27] {Bjorck1996}}
& \twolines {\hspace*{0.25em} \| \DA \|_2 \hspace* {0.25em}} {\scaleA = 1}
& \displaystyle \Db = 0
& \twolines {\| \Dr \|_2} {\scaleR = 1}
& \displaystyle {\| r \|_2 \over \sigmamin} + \| x \|_2
& 2
\\ \hline
\leftbox {\citet [p.\ 655] {Stewart1977b}, \citet [p.\ 160, sec.\ 5.2] {Stewart1990}}
& \twolines {\| \DA \|_2} {\scaleA = 1}
& \displaystyle \Db = 0
& \twolines {\| \Dr \|_2} {\scaleR = 1}
& \displaystyle {\| b \|_2 \over \sigmamin} 
&\sqrt 2 \, \boldkappa
\\ \hline
\leftbox {\citet [p.\ 242, eqn.\ 5.3.9] {Golub1996}, \citet [p.\ 382, eqn.\ 20.2] {Higham2002}$^{**}$}
& \multicolumn {2} {c |} {\twolines {\hspace*{-0.66em} \max \left\{ {\| \DA \|_2 \over \| A \|_2}, {\|
\Db \|_2 \over \| b \|_2} \right\} \hspace*{-0.66em}} {\scaleA = \| A \|_2 \quad \scaleB = \| b \|_2}}
& \twolines {\fracstrut {\| \Dr \|_2 \over \| b \|_2} \fracstrut} {\scaleR = \| b \|_2}
& \displaystyle 2 \, {\| A \|_2 \over \sigmamin} + 1^{**}
& \boldkappa
\\ \hline
\leftbox {Theorem \ref {thm:condition-numbers},\\ equation (\ref {eqn:thm})\vrule depth3ex height0pt width\strutwidth}
& \displaystyle {\| \DA \|_2 \over \vphantom {(} \scaleA}
& \Db = 0
& \displaystyle {\| \Dr \|_2 \over \scaleR}
& \vrule depth3.5ex height5.5ex width\strutwidth \displaystyle {{\scaleA \over \scaleR} \, \sqrt {\left({\| r \|_2 \over \sigmamin}\right)^2 + \| x \|_2^2}}
&\sqrt 2
\\ \hline
\multicolumn {6} {l} {\footnotesize \hspace*{1em} \begin {minipage} {4.75in} \raggedright \textstrut \llap {$\;^{**}\;$}The formula of Golub and van Loan and of Higham amounts to an estimate for $\chi_r (A) + \chi_r (b)$\\[-0.5ex] and is compared against the sum of $\chi_r (b)$ and the tight estimate for $\chi_r (A)$. See section \ref {sec:GVLH}.\end {minipage}}
\end {array}
\end {displaymath}
\end {table}

\subsection {The estimate of Stewart} 

The value of \citet [p.\ 655] {Stewart1977b} is also reported by \citet [p.\ 160, sec.\ 5.2] {Stewart1990}. It again is for choices $\Db = 0$ and $\scaleA = \scaleR = 1$. Some assembly is required. Let $B = A + \DA$ be the perturbed matrix. Assume $\| \DA \|_2 < \sigmamin$ so that $B$ also has full rank. 

For any matrix $M$, let $P_M = M M^\dag$ be the orthogonal projection into the column space of $M$. With $\Db = 0$ the difference between the residuals of the original and the perturbed problems (\ref {eqn:original-problem}, \ref {eqn:perturbed-problem}) is $\Dr = (I - P_B) b - (I - P_A) b$ so it is always true that
\begin {equation}
\label {eqn:stewart-1}
\| \Dr \|_2 \le \| P_A - P_B \|_2 \, \| b \|_2 \, .
\end {equation}
\citet [p.\ 655] {Stewart1977b} remarks that $\| P_A - P_B \|_2$ is to be bounded by applying an earlier result. He does not intend $\| P_A - P_B \|_2 < 1$ (p.\ 651, eqn.\ 4.1) which converts equation (\ref {eqn:stewart-1}) into the useless $\| \Dr \|_2 < \| b \|_2$. Stewart means a complicated expression that introduces $\| \DA \|_2$ into the bound. This expression requires some preparation that is more easily followed in the presentation of \citet [pp.\ 160, 153, 148, 137] {Stewart1990}.

Continuing the assembly of the bound, let the singular value decomposition of $A$ be
\begin {displaymath}
{\onetwo {U_1} {U_2}}^t A V = \twoone {A_1} 0
\end {displaymath}
where $[\,{U_1},\, {U_2}]$ and $V$ are square orthonormal matrices and where $A_1$ is the square diagonal  matrix of singular values. Let the corresponding factorization of $\DA$ be \citep [p.\ 137] {Stewart1990}
\begin {displaymath}
{\onetwo {U_1} {U_2}}^t \DA V = \twoone {E_1} {E_2}
\end {displaymath}
where $\| E_i \|_2 \le \| \DA \|_2$ for $i = 1,2$. \citet [p.\ 148] {Stewart1990} define
\begin {displaymath}
\hat \kappa = \| A \|_2 \, \| (A_1 + E_1)^{-1} \|_2 \, .
\end {displaymath}
From the triangle inequality and from the Neumann series expansion for $(A_1 + E_1)^{-1}$,
\begin {displaymath}
\left| \, \| A_1^{-1} \|_2 - \| (A_1 + E_1)^{-1} \|_2 \, \right| \le \| A_1^{-1} - (A_1 + E_1)^{-1} \|_2 \le {\mathcal O} (\| \DA \|_2^2) \, .
\end {displaymath}
These last two equations combine to 
\begin {equation}
\label {eqn:stewart-2}
\hat \kappa = \| A \|_2 \, \sigmamin^{-1} + {\mathcal O} (\| \DA \|_2^2 ) \, .
\end {equation}

The final step applies a bound that requires some further hypotheses. For any matrix $M$, similar to $P_M = M M^\dag$, let $R_M = P_{M^t} = (M^\dag M)^t$ be the orthogonal projection into the row space of $M$ (viewing the rows as column vectors). Since $B^\dag = (B^t B)^{-1} B^t$ is a continuous function of $\DA$, therefore both $B B^\dag - A A^\dag$ and $B^\dag B - A^\dag A$ converge to $0$ as $\DA$ approaches $0$.\footnote {Specific bounds on the norms of these differences can be derived from \citet [p.\ 221, thm.\ 4.1] {Wedin1973}.}  If $\| \DA \|_2$ is sufficiently small that both $\| P_A - P_B \|_2 \le 1$ and $\| R_A - R_B \|_2 \le 1$, then it can be shown \citep [p.\ 153, eqn.\ 4.1] {Stewart1990}
\begin {equation}
\label {eqn:stewart-3}
\| P_A - P_B \|_2  \le {\hat \kappa \, \| E_2 \|_2 / \| A \|_2 \over \left[ 1 + (\hat \kappa \, \| E_2 \|_2 / \| A \|_2)^2 \right]^{1/2}} \, .
\end {equation}
Altogether, combining equations (\ref {eqn:stewart-1}--\ref {eqn:stewart-3}) leaves
\begin {equation}
\label {eqn:stewart-4}
\| \Dr \|_2 \le {\| b \|_2 \over \sigmamin} \, \| \DA \|_2 + {\mathcal O} (\| \DA \|_2^2) \, ,
\end {equation}
which is the bound from which the condition estimate in table \ref {tab:published} is taken.

The estimate for $\chi_r (A)$ in equation (\ref {eqn:stewart-4}) can be obtained from theorem \ref {thm:condition-numbers} by increasing the second term in the upper bound (\ref {eqn:thm}) by the factor $\vds^2$,
\begin {displaymath}
\sqrt {\left({\| r \|_2 \over \sigmamin}\right)^2 + \| x \|_2^2}
\le
\sqrt {\left({\| r \|_2 \over \sigmamin}\right)^2 + \vds^2 \| x \|_2^2}
=
\sqrt {{\| r \|_2^2 + \| A x \|_2^2 \over \sigmamin^2}}
=
{\| b \|_2 \over \sigmamin} \, .
\end {displaymath}
Consequently, Stewart and Sun's value can overestimate the upper bound for $\chi_r (A)$ by as much as $\vds$ depending on circumstances. The worst situation is illustrated by the example of section \ref {sec:example} with $\phi = 0$ and $\beta \gg 1 / \alpha$,
\begin {displaymath}
{\displaystyle {\| b \|_2 \over \strut \sigmamin} \over \displaystyle \sqrt {\left({\| r \|_2 \over \sigmamin}\right)^2 + \| x \|_2^2}}
=
{\displaystyle {\sqrt {1 + \beta^2} \over \strut \alpha} \over \displaystyle \sqrt {
{1 \over \alpha^2} + \beta^2 \cos^2 (\phi) + {\beta^2 \over \alpha^2} \sin^2 (\phi)
}}
=
{\sqrt {1 + \beta^2} \over \sqrt {1 + \alpha^2 \beta^2}}
\approx
{1 \over \alpha}
=
\boldkappa \, .
\end {displaymath}

\subsection {The estimate of Golub and Van Loan and of Higham}
\label {sec:GVLH}

\newcommand {\maxformula} {\max \left\{  {\| \DA \|_2 \over \| A \|_2}, \, {\| \Db \|_2 \over \| b \|_2} \right\}}

The condition estimate of \citet [p.\ 242, eqn.\ 5.3.9] {Golub1996} and of \citet [p.\ 382, eqn.\ 20.2] {Higham2002} is for the choices $\DA \ne 0$ and $\Db \ne 0$ with the scale factors $\scaleA = \| A \|_2$ and  $\scaleR = \scaleB = \| b \|_2$. They take the approach of equation (\ref {eqn:single}) that uses a single quantity, $\varepsilon$, to measure the perturbations to $A$ and $b$,
\begin {equation}
\label {eqn:epsilon}
\left. \begin {array} {r c l} \| \DA \|_2& \le& \varepsilon \| A \|_2\\ \noalign {\smallskip} \| \Db \|_2& \le& \varepsilon \| b \|_2\end {array} \right\} \quad \mbox {equivalently} \quad
\varepsilon = \maxformula . 
\end {equation}
Since $\scaleB = \scaleR$, this approach can transform the bound (\ref {eqn:differential-bound}) as follows,
\begin {eqnarray}
\nonumber
{\| \Dr \|_2 \over \scaleB}& \le& \| \JrA \| \, {\scaleA \over \scaleB} \, {\| \DA \|_2 \over \scaleA} + \| J_r (b) \| \, {\| \Db \| \over \scaleB} + o ( \varepsilon )\\
\noalign {\smallskip}
\nonumber
\label {eqn:GVLH-1}
& \le & \left[  \| \JrA \| \, {\scaleA \over \scaleB} + \| J_r (b) \| \right] \varepsilon + o  (\varepsilon)\\
\noalign {\smallskip}
\label {eqn:GVLH-2}
& = & \big[  \chi_r (A) + \chi_r (b) \big] \varepsilon + o (\varepsilon) \, .
\end {eqnarray}
From theorem \ref {thm:condition-numbers} with the choices $\scaleR = \scaleB = \| b \|_2$,
\begin {eqnarray}
\nonumber
\chi_r (A) + \chi_r (b)
& \le& \left\{ {\| A \|_2 \over \| b \|_2} \sqrt {\left({\| r \|_2 \over \sigmamin}\right)^2 + \| x \|_2^2} \right\} + 1\\
\noalign {\smallskip}
\label {eqn:GVLH-3}
& \le& \left\{ {\| A \|_2 \over \| b \|_2} \sqrt {\left({\| r \|_2 \over \sigmamin}\right)^2 + \vds^2 \| x \|_2^2} \right\} + 1  = \boldkappa + 1 \, .
\end {eqnarray}

Golub and Van Loan and Higham state a larger value, $2 \boldkappa + 1$.\footnote {Their value resembles the $\sqrt 2 \, \boldkappa + 1$ that was originally stated by \citet [p.\ 16, eqn.\ 7.7] {Bjorck1967a}.} Since these formulas can be derived from the sum $\chi_r (A) + \chi_r (b)$ they are not joint condition numbers in the sense of equation (\ref {eqn:single}).
Moreover, the derivation inserts $\vds$ into equation (\ref {eqn:GVLH-3}), so the result can overestimate the sum by as much as a factor of $\boldkappa$. Close to the worst situation for the specific value $2 \boldkappa + 1$ of Golub, Van Loan and Higham is again illustrated by the example of section \ref {sec:example} with $\phi = 0$ and $\beta \gg 1 / \alpha$,
\begin {displaymath}
{\displaystyle {\strut 2 \boldkappa + 1} \over \displaystyle \left\{ {\| A \|_2 \over \| b \|_2} \sqrt {\left({\| r \|_2 \over \sigmamin}\right)^2 + \| x \|_2^2} \right\} + 1}
=
{\displaystyle {{2 \over \strut \alpha} + 1} \over \displaystyle \left\{ {1 \over \sqrt {1 + \beta^2}} \sqrt {
{1 \over \alpha^2} + \beta^2} \right\} + 1}
\approx
{1 \over \alpha}
=
\boldkappa \, .
\end {displaymath}

Note the bound (\ref {eqn:GVLH-3}) is not sensitive to the angle $\boldtheta$ between $r$ and $\col (A)$ because of the choice for the scale factor $\scaleR = \| b \|_2$. Choices for $\scaleR$ are discussed in section \ref {sec:residual-scaling}.

\section {Discussion}
\label {sec:discuss}

\subsection {Measuring perturbations to $\mathbf r$ relative to $\mathbf b$}
\label {sec:residual-scaling}

As mentioned in section \ref {sec:norms}, scaled changes to the residual are typically measured by choosing $\scaleR = \| r \|_2$ or $\| b \|_2$. The two cases are contrasted in Table \ref {tab:scaling}. The choice $\scaleR = \| b \|_2$ makes it appear that $\boldtheta$ is not a source of ill-conditioning because the sensitivity of $r$ to $A$ is masked by measuring changes to $r$ against the always larger vector $b$. The choice $\scaleR = \| r \|_2$ measures perturbations relatively. The next section \ref {sec:iterative} describes a situation when the relative measure is appropriate.

\begin {table} [h!]
\caption {\it Effect of scaling on condition numbers for the least squares residual. The full rank least squares problem is $\min_x \| b - A x \|_2$, the solution is $x$, the residual is $r$, $\boldkappa = \| A \|_2 / \sigmamin$ is the spectral matrix condition number of $A$, $\sigmamin$ is the smallest singular value of $A$, $\vds = \| A x \|_2 / (\| x \|_2 \, \sigmamin)$ is van der Sluis's ratio between $1$ and $\boldkappa$, and $\boldtheta$ is the angle between $b$ and $\col (A)$.}
\label {tab:scaling}
\newcommand {\textstrut} {\vrule depth1.125ex height2.375ex width\strutwidth}
\newcommand {\topstrut} {\vrule depth0ex height2.75ex width\strutwidth}
\newcommand {\botstrut} {\vrule depth1.0ex height0ex width\strutwidth}
\newcommand {\equationstrut} {\vrule depth3.5ex height5.5ex width\strutwidth}
\newcommand {\leftbox} [1] {\begin {minipage}{8.5em}\footnotesize \vrule depth0ex height3.0ex width\strutwidth \raggedright #1\vrule depth1.5ex height2ex width\strutwidth \end {minipage}}
\setlength {\arraycolsep} {0.5em}
\vspace {-2ex}
\begin {displaymath}
\small
\begin {array} {| c | c | c | c | c |}
\hline
\multicolumn {3} {| c |} {\textstrut \mbox {norms and scale factors}}
& \multicolumn {2} {c |} {\textstrut \mbox {condition numbers}}
\\
\cline {1-5} 
\displaystyle {\topstrut \| \DA \|_2 \over \botstrut \scaleA}
& \displaystyle {\topstrut \| \Db \|_2 \over \botstrut \scaleB}
& \displaystyle {\topstrut \| \Dr \|_2 \over \botstrut \scaleR}
& \mbox {tight estimate for $\chi_r (A)$}
& \chi_r (b)
\\ \hline \hline
\scaleA = \| A \|_2
& \scaleB = \| b \|_2
& \scaleR = \| r \|_2
& \equationstrut \displaystyle \boldkappa \sqrt { 1 + \left( \cot (\boldtheta) \over \vds \right)^2}
& \csc (\boldtheta)
\\
\scaleA = \| A \|_2
& \scaleB = \| b \|_2
& \scaleR = \| b \|_2
& \equationstrut \displaystyle \boldkappa \sqrt { \sin^2 (\boldtheta) + \left( \cos (\boldtheta) \over \vds \right)^2}
& 1
\\ \hline
\end {array}
\end {displaymath}
\end {table}

\subsection {Significance for iterative methods}
\label {sec:iterative}

Many iterative methods proceed by building orthogonal bases from the residuals of least squares projections. For example, for a symmetric matrix $A$ and a unit vector $v_1$, the Lanczos iteration
\begin {displaymath}
\beta_{j+1} v_{j+1} = A v_j - \alpha_j v_j - \beta_j v_{j-1}
\end {displaymath}
produces a sequence of orthonormal vectors $v_1$, $v_2$, $v_3$, \dots. This algorithm can be viewed as repeatedly evaluating a residual $r_{j+1} = \beta_{j+1} v_{j+1}$ for either of two orthogonal projections: (1) the projection of $A v_j$ into the span of $v_{j-1}$ and $v_j$, or (2) the orthogonal projection of $A^j v_1$ into the Krylov subspace spanned by $v_1$, $A v_1$, \dots, $A^{j-1} v_1$. 

The appropriate scale factor for measuring perturbations to $r_{j+1}$ is $\scaleR = \| r_{j+1} \|_2$. The relative error in $r_{j+1}$ becomes the absolute error in the normalized vector, $v_{j+1}$, which continues the Lanczos iteration. In ideal circumstances the vectors $v_{j-1}$ and $v_j$ are close to orthonormal. If $A = [ v_{j-1}, v_j]$ is an orthonormal matrix, then $\boldkappa = \vds = 1$ so the tight estimate in Table \ref {tab:scaling} simplifies to $\chi_r (A) \approx \csc (\boldtheta) = \chi_r (b)$ where $\boldtheta$ is the angle between $b = A v_j$ and $\col (A)$. Thus $r_{j+1}$ is ill-conditioned when $\boldtheta$ is small.

\subsection {Condition numbers of orthogonal projections} 

In the least squares problem, the condition number of the orthogonal projection $Ax$ is essentially that of the residual. In addition to equations (\ref {eqn:norms}, \ref {eqn:scale-factors}), it is necessary to specify the scaled norm for the projection:
\begin {displaymath}
\| \DAx \| = {\| \DAx \|_2 \over \scaleP} \quad \mbox {where} \quad \scaleP = \| Ax \|_2 \, .
\end {displaymath}
From $Ax = P b$ follows $J_{Ax} (b) = P$ hence $\chi_{Ax} (b) = \scaleB / \scaleP$. Since $Ax = b - r$ so $J_{Ax} (A) = - J_r (A)$. With $\scaleP$ replacing $\scaleR$ in the formulas, the condition numbers of the orthogonal projection are
\begin {displaymath}
\setlength {\arraycolsep} {0.33em}
\begin {array} {r c c c l}
&& \displaystyle \chi_{Ax} (b)
& =
& \sec (\boldtheta) \, ,
\\ \noalign {\bigskip}
\displaystyle {{1 \over \sqrt 2} \, {\boldkappa \sqrt { \tan^2 (\boldtheta) + {1 \over \vds^2}}}}
& \le 
& \displaystyle \chi_{Ax} (A)
& \le 
& \displaystyle {\boldkappa \sqrt { \tan^2 (\boldtheta) + {1 \over \vds^2}}} \, .
\end {array}
\end {displaymath}
Both $\boldkappa$ and $\boldtheta$ are independent sources of ill-conditioning. 

\subsection {Column transformations} 

The linear least squares residual is invariant with respect to transformations of the matrix columns, so there is reason to seek changes to the columns that might reduce the condition number of the residual. If $A$ is replaced by $A M$ for some nonsingular matrix $M$ that makes $AM$ an orthonormal matrix, then with the scale factors of equation (\ref {eqn:scale-factors})  it has been noted in section \ref {sec:iterative} that the tight estimate is $\chi_r (AM) \approx \csc (\boldtheta)$ which leaves only $\boldtheta$ as a source of ill-conditioning. 

A less costly transformation is $M = D$ for a diagonal matrix $D$. Two reasons suggest choosing $D$ to equilibrate the columns of $AD$. First, least squares problems typically are solved using the $QR$ factorization. The errors of that calculation can be accounted for by backward rounding errors whose relative size in each column is roughly the same across all columns \citep [p.\ 385, thm.\ 20.3] {Higham2002}. Second, equilibrating the columns is approximately the optimal column scaling to reduce the matrix condition number \citep {vanderSluis1969}. Nevertheless, even if $\boldkappa (AD) \le \boldkappa (A)$, the scaling also alters van der Sluis's ratio in equation (\ref {eqn:three}), so it is unclear whether the net change to the condition number in equation (\ref {eqn:chiA}) is for the better.

\section* {Acknowledgements}

I thank the editor Prof.\ D.\ O'Leary and the three referees for corrections and suggestions that much improved this paper.




\frenchspacing
\raggedright

\end {document}